\documentclass[11pt,a4paper,reqno]{amsart}
\usepackage[hmargin=28mm,top=21mm,bottom=22mm]{geometry}
\usepackage[T1]{fontenc}
\usepackage[utf8]{inputenc}
\usepackage{lmodern}
\usepackage{amssymb,enumitem}
\usepackage{xcolor}
\usepackage[unicode,colorlinks=true,linkcolor=black,citecolor=black,urlcolor=blue]{hyperref}
\hypersetup{pdftitle={Diagonal L2 Dolbeault cohomology of pseudoconvex domains in complex projective space},pdfauthor={Masakazu Takakura}}
\newtheorem{theorem}{Theorem}[section]
\newtheorem{lemma}[theorem]{Lemma}
\newtheorem{corollary}[theorem]{Corollary}
\newtheoremstyle{examplebreak}{}{}{\normalfont}{}{\bfseries}{.}{\newline}{}
\theoremstyle{examplebreak}
\newtheorem{example}[theorem]{Example}
\makeatletter
\renewcommand{\subsection}{\@startsection{subsection}{2}{\z@}%
  {.65\linespacing\@plus.2\linespacing}{.25\linespacing}%
  {\normalfont\bfseries}}
\makeatother
\newcommand{\condref}[1]{\textnormal{\ref*{cond:#1}}}
\newcommand{\plaincite}[2][]{\begin{NoHyper}\cite[#1]{#2}\end{NoHyper}}
\newcommand{\PP}{\mathbb P}
\newcommand{\CC}{\mathbb C}
\newcommand{\ddbar}{\bar\partial}
\newcommand{\FS}{\mathrm{FS}}
\newcommand{\Dom}{\operatorname{Dom}}

\newcommand{\Vol}{\operatorname{Vol}}
\newcommand{\hL}[2]{H^{#1,#2}_{\smash[t]{\scriptscriptstyle(2)}}}
\setlist[enumerate]{leftmargin=2.7em,topsep=.4em,itemsep=.1em,parsep=0pt}
\title[Diagonal $L^2$ Dolbeault cohomology]{Diagonal $L^2$ Dolbeault cohomology of pseudoconvex domains in complex projective space}
\author{Masakazu Takakura}
\date{}
\begin{document}
\raggedbottom

\begin{abstract}
Let $\Omega$ be a proper pseudoconvex domain in complex projective space $\PP^n$, equipped with the restriction of the Fubini--Study metric. We prove that the following three conditions are equivalent: the $L^2$ Dolbeault cohomology of $\Omega$ vanishes in every positive diagonal bidegree, it vanishes in some positive diagonal bidegree, and the complement of $\Omega$ is not locally polar.
\end{abstract}
\maketitle

\section{Introduction}

On a Hermitian manifold $(Y,g)$, let $L^2_{p,q}=L^2_{p,q}(Y,g)$ denote the space of square-integrable $(p,q)$-forms. We take $\ddbar_{p,q}$ to be the distributional $\ddbar$ operator with domain
\[
 \Dom\ddbar_{p,q}=\{u\in L^2_{p,q}:\ddbar u\in L^2_{p,q+1}\}.
\]
For $q>0$, the $L^2$ Dolbeault cohomology group is
\[
 \hL pq(Y,g)
 =\frac{\ker\bigl(\ddbar_{p,q}:\Dom\ddbar_{p,q}\longrightarrow L^2_{p,q+1}\bigr)}
        {\ddbar_{p,q-1}(\Dom\ddbar_{p,q-1})}.
\]

The $L^2$ Dolbeault cohomology of pseudoconvex domains in complex projective space $\PP^n$ has been studied by several authors. Let $\Omega$ be a proper pseudoconvex domain in $\PP^n$. We write $\hL pq(\Omega)$ for $\hL pq(\Omega,g_{\FS}|_\Omega)$, where $g_{\FS}$ is the Fubini--Study metric with Kähler form $\omega=i\partial\ddbar\log(1+|z|^2)$ in affine coordinates.

For $n\ge2$, Henkin and Iordan \plaincite[Corollary~3.6, p.~861]{HI} proved that if $\partial\Omega$ is of class $C^2$, then $\hL pq(\Omega)=0$ for every bidegree $(p,q)$ with $q>0$. The same vanishing holds for Lipschitz boundaries by Fu and Shaw \plaincite[Lemma~4.1]{FS22}. We are interested in the behavior of these groups without any regularity assumption on $\partial\Omega$.
\begin{enumerate}[label=\textup{(\roman*)}]
\item For $p<q$, the vanishing $\hL pq(\Omega)=0$ is a direct consequence of Demailly's $L^2$ existence theorem and the Fubini--Study curvature calculation in Lemma~\ref*{lem:curv} (see Corollary~\ref*{cor:p-less-q}). No boundary regularity is required.
\item If $p>q$, vanishing can fail. Shaw's example is the projective Hartogs triangle
\[
 \mathbb H=\{[Z_0:Z_1:Z_2]\in\PP^2: |Z_1|<|Z_2|\},
\]
whose boundary is not Lipschitz at $[1:0:0]$. Shaw showed that $\hL 21(\mathbb H)$ is infinite-dimensional and non-Hausdorff \plaincite[Lemma~4.3(3) and Corollary~5.2, pp.~15--16]{Shaw}.
\item The diagonal $p=q$ is the borderline case. Shaw asked whether $\hL 11(\mathbb H)$ is Hausdorff and whether it vanishes \plaincite[\S6, Problem~5, p.~19]{Shaw}. Recently, J.~Chen answered both questions in the affirmative by proving that, for every $n\ge2$, the projective Hartogs triangle
\[
 \mathbb H_n=\{[Z_0:\cdots:Z_n]\in\PP^n:|Z_{n-1}|<|Z_n|\}
\]
satisfies $\hL 11(\mathbb H_n)=0$ \plaincite[(1.2), p.~1, and Corollary~1.4, p.~3]{Chen}.
\end{enumerate}

More generally, given a proper pseudoconvex domain $\Omega\subsetneq\PP^n$ and $1\le p\le n$, when does $\hL pp(\Omega)$ vanish?~A basic example of diagonal nonvanishing occurs when the complement $K=\PP^n\setminus\Omega$ is locally polar. Recall that $K$ is \emph{locally polar} (respectively, \emph{locally pluripolar}) if each point of $K$ has a connected coordinate neighborhood $U$ and a subharmonic (respectively, plurisubharmonic) function $v\not\equiv-\infty$ on $U$ such that $K\cap U\subset\{v=-\infty\}$. Subharmonicity here is taken in real Euclidean coordinates. If $K$ is locally polar, applying B.-Y.~Chen, Wu and Wang \plaincite[Corollary~2.4]{CWW} in coordinate charts shows that the distributional equation $\ddbar u=v$ extends across $K$ for $L^2$ forms $u$ and $v$. Since $K$ has measure zero, restriction thus identifies the $L^2$ Dolbeault complexes of $\PP^n$ and $\Omega$ isometrically. In particular,
\[
 \hL pp(\Omega)=\CC[\omega^p|_\Omega],\qquad 1\le p\le n .
\]
Our main result shows that this is the only obstruction to diagonal vanishing.

\begin{samepage}
\begin{theorem}\label{thm:main}
Let $\Omega\subsetneq\PP^n$ be a pseudoconvex domain and put $K=\PP^n\setminus\Omega$. The following conditions are equivalent.
\begin{enumerate}[label=\textup{(\roman*)}]
 \item\label{cond:nonpolar} $K$ is not locally polar.
 \item\label{cond:nonpluripolar} $K$ is not locally pluripolar.
 \item\label{cond:some} $\hL pp(\Omega)=0$ for some $1\le p\le n$.
 \item\label{cond:all} $\hL pp(\Omega)=0$ for every $1\le p\le n$.
\end{enumerate}
\end{theorem}
\end{samepage}

The equivalence of \condref{nonpolar} and \condref{nonpluripolar} follows from Sadullaev's theorem that a polar pseudoconcave set is pluripolar \plaincite[Theorem~2, p.~92]{Sad}, applied in affine charts. To prove $\condref{nonpluripolar}\Rightarrow\condref{all}$, we combine Demailly's $L^2$ existence theorem with the Poisson equation on $\PP^n$.

\subsection*{Acknowledgements}
{\raggedright This work was supported by JST SPRING, Grant Number \mbox{JPMJSP2156}.\par}

\subsection*{Use of AI}
In preparing this paper, the author used ChatGPT (OpenAI, GPT-6.2 Astra) to edit the English text and search the literature. For the construction of $\tau$ in the proof of Lemma~\ref*{lem:mean}, ChatGPT suggested using a square-integrable holomorphic function on $B\cap\Omega$ that does not extend to $B$ (see the slicing argument in Lemma~\ref*{lem:berg}). All other mathematical content is due to the author.

\section{Preliminaries}

\subsection{Demailly's \texorpdfstring{$L^2$}{L2} existence theorem}

Let $(Y,\omega)$ be a Kähler manifold of dimension $n$ and $(E,h)$ a holomorphic Hermitian vector bundle. Denote by $\Theta(E,h)=D_h^2$ the curvature of the Chern connection $D_h$ on $(E,h)$. Write $dV_\omega=\omega^n/n!$, and let $\Lambda_\omega$ denote the adjoint of wedging with $\omega$. On $E$-valued $(n,q)$-forms the curvature operator is
\[
 B_q=[i\Theta(E,h),\Lambda_\omega]
     =(i\Theta(E,h)\wedge\cdot)\,\Lambda_\omega .
\]
We write $B_q\ge0$ if
\[
 \langle B_q(x)\xi,\xi\rangle_{h,\omega}\ge0
 \qquad\text{for all }x\in Y
 \text{ and }\xi\in\Lambda^{n,q}T_x^*Y\otimes E_x .
\]
We use
\plaincite[Proposition~3.12]{Dem16}
with $\eta=1$ and $\lambda=\varepsilon$, passing to a weak limit as $\varepsilon\downarrow0$.

\begin{samepage}
\begin{theorem}[Demailly]\label{thm:demailly}
Suppose $Y$ admits a complete Kähler metric. Fix $1\le q\le n$ and assume that $B_q\ge0$ for the given Kähler metric $\omega$. If $g\in L^2_{n,q}(Y,E;\omega,h)$ is $\ddbar$-closed and
\[
 M(g):=\lim_{\varepsilon\downarrow0}\int_Y
 \big\langle(B_q+\varepsilon\operatorname{Id})^{-1}g,g\big\rangle_{h,\omega}
 \,dV_\omega <\infty ,
\]
there exists $u\in L^2_{n,q-1}(Y,E;\omega,h)$ such that $\ddbar u=g$ and $\|u\|_{L^2(\omega,h)}^2\le M(g)$.
\end{theorem}
\end{samepage}

\subsection{The Fubini--Study curvature}
Put $E_p=\bigwedge^{n-p}T_{\PP^n}$ with the metric induced by $\omega$. For an $E_{p,x}$-valued $(n,q)$-form $u$ at $x$, let $\widehat u$ denote the ordinary $(p,q)$-form obtained by contraction under the isometric isomorphism $K_{\PP^n}\otimes E_p\simeq\Omega^p_{\PP^n}$. Thus, pointwise,
\[ u\in E_{p,x}\otimes\Lambda^{n,q}T_x^*\PP^n\ \longmapsto\ \widehat u\in\Lambda^{p,q}T_x^*\PP^n .\]
Write $B_{p,q}\widehat u=\widehat{B_qu}$. On ordinary forms set $L=\omega\wedge\cdot$ and $\Lambda=L^*$.

\begin{lemma}\label{lem:curv}
We have $B_{p,q}=q(n-p+1)\operatorname{Id}-L\Lambda$. In particular, if
$0\le p\le q\le n$, then
\[
 B_{p,q}\ge(n+1)(q-p)\operatorname{Id}.
\]
For $1\le p\le n$, the pointwise kernel of $B_{p,p}$ is $\CC\omega^p$, and
\[
 (\ker B_{p,p})^\perp
 =\ker\bigl(\alpha\longmapsto\alpha\wedge\omega^{n-p}\bigr)
 \quad\bigl(\alpha\in\Lambda^{p,p}T_{\PP^n}^*\bigr).
\]
On this orthogonal complement, $B_{p,p}\ge n\operatorname{Id}$.
\end{lemma}
\begin{proof}
Put $T=T_x^{1,0}\PP^n$ and $m=n-p$. Choose a unitary basis $(e_a)$ of $T$ with dual basis $(\theta^a)$, and write $i\Theta(E_p)_x=i\sum_{a,b}\theta^a\wedge\bar\theta^b\otimes R_{a\bar b}$. By the calculation of Henkin and Iordan \plaincite[proof of Lemma~3.3, p.~860, following (3.2)]{HI}, the Fubini--Study curvature acts on a basis of $\bigwedge^m T$ by
\[
 R_{a\bar b}e_I=m\delta_{ab}e_I+S_{ab}e_I,
 \qquad e_I=e_{i_1}\wedge\cdots\wedge e_{i_m}\quad(i_1<\cdots<i_m).
\]
Here $\delta_{ab}$ is $1$ if $a=b$ and $0$ otherwise. On a basis wedge of any degree, $S_{ab}$ replaces $e_b$ by $e_a$ in its original position. The result is zero if $e_b$ is absent or the replacement gives a repeated factor.

Thus $i\Theta(E_p)_x=A+B$, where
\[
 A=m\omega\otimes\operatorname{Id}_{E_{p,x}},\qquad
 B=i\sum_{a,b}\theta^a\wedge\bar\theta^b\otimes S_{ab}.
\]
We identify $A$ and $B$ with their wedge operators and compute $[A,\Lambda]$ and $[B,\Lambda]$ separately on $E_{p,x}$-valued $(n,q)$-forms. Since $u$ has holomorphic degree $n$, $Au=Bu=0$.

First, $A$ acts by $mL$. Since $[L,\Lambda]=(k-n)\operatorname{Id}$ in total degree $k$, we have
\[
 [A,\Lambda]u=m[L,\Lambda]u=mq\,u.
\]

For $[B,\Lambda]$, it suffices by linearity to take $u=v\otimes(\tau\wedge\gamma)$, where $v\in\bigwedge^m T$, $\tau=\theta^1\wedge\cdots\wedge\theta^n$, and $\gamma$ is a $(0,q)$-form. Expanding $B\Lambda$ in the unitary coframe gives (cf.\ \plaincite[(3.6), p.~469]{Dem82})
\[
 [B,\Lambda]u=B\Lambda u
 =\sum_{a,b}S_{ab}v\otimes
   (\tau\wedge\bar\theta^b\wedge\iota_{\bar e_a}\gamma).
\]
Put $\beta=\iota_v\tau$, so that $\widehat u=\beta\wedge\gamma$. For $w\in\bigwedge^p T$, the derivation $S_{ab}$ acts on $v\wedge w\in\det T$ by $\delta_{ab}$. Evaluating $\tau$ on $(S_{ab}v)\wedge w=\delta_{ab}\,v\wedge w-v\wedge S_{ab}w$ gives
\[
 \iota_{S_{ab}v}\tau=\delta_{ab}\beta-\theta^b\wedge\iota_{e_a}\beta.
\]
Taking hats in the expression for $[B,\Lambda]u$, we obtain
\[
 \begin{aligned}
 \widehat{[B,\Lambda]u}
 &=\sum_{a,b}(\delta_{ab}\beta-\theta^b\wedge\iota_{e_a}\beta)
       \wedge(\bar\theta^b\wedge\iota_{\bar e_a}\gamma)\\
 &=q\,\widehat u-L\Lambda\widehat u .
 \end{aligned}
\]
Here $\sum_a\bar\theta^a\wedge\iota_{\bar e_a}\gamma=q\gamma$. In the remaining sum the contractions over $a$ give $\Lambda$ and the wedges over $b$ give $L$, since $L=i\sum_b\theta^b\wedge\bar\theta^b\wedge\cdot$ and $\Lambda=-i\sum_a\iota_{\bar e_a}\iota_{e_a}$. Adding the two commutators gives
\[
 \widehat{B_qu}=\widehat{[A,\Lambda]u}+\widehat{[B,\Lambda]u}
     =q(n-p+1)\widehat u-L\Lambda\widehat u .
\]

Finally write $\widehat u=\sum_r L^r\xi_r$ for its orthogonal Lefschetz decomposition, where $\xi_r$ is a primitive $(p-r,q-r)$-form and $\max(0,p+q-n)\le r\le\min(p,q)$. Since $\Lambda\xi_r=0$, the commutator $[\Lambda,L]=(n-k)\operatorname{Id}$ on forms of total degree $k$ gives $L\Lambda(L^r\xi_r)=r(n-p-q+r+1)L^r\xi_r$. Thus
\[
 \begin{aligned}
 |u|^2=|\widehat u|^2&=\sum_r|L^r\xi_r|^2,\\
 \langle B_qu,u\rangle=\langle B_{p,q}\widehat u,\widehat u\rangle
  &=\sum_r (q-r)(n+1-p+r)\,|L^r\xi_r|^2 .
 \end{aligned}
\]
For $p\le q$, the coefficient exceeds $(n+1)(q-p)$ by $(p-r)(n+1-q+r)\ge0$, with equality exactly at $r=p$. Thus equality in the curvature bound holds precisely when $\widehat u=L^p\xi_p$. For $q=p$, the zero eigenspace is $L^p\CC=\CC\omega^p$. The other eigenvalues are $s(n+1-s)\ge n$, with $s=p-r\ge1$. Finally, $*\omega^p=\frac{p!}{(n-p)!}\omega^{n-p}$ identifies $(\CC\omega^p)^\perp$ with the kernel of wedging with $\omega^{n-p}$.
\end{proof}

The positive lower bound in Lemma~\ref*{lem:curv} gives the following vanishing theorem for domains admitting a complete Kähler metric.
\begin{corollary}\label{cor:p-less-q}
Let $D\subset\PP^n$ be a domain admitting a complete Kähler metric. With respect to the induced Fubini--Study metric, $\hL pq(D)=0$ for $0\le p<q\le n$.
\end{corollary}
\begin{proof}
Apply Theorem~\ref*{thm:demailly} to $E_p|_D$ with the lower bound in Lemma~\ref*{lem:curv}.
\end{proof}
In the borderline case $p=q$, one still has the following solvability result.
\begin{corollary}\label{cor:orth}
If $1\le p\le n$ and $g\in L^2_{p,p}(\Omega)$ satisfies $\ddbar g=0$ and $g\wedge\omega^{n-p}=0$, there exists $b\in L^2_{p,p-1}(\Omega)$ with $\ddbar b=g$ and $\|b\|_{L^2}^2\le n^{-1}\|g\|_{L^2}^2$.
\end{corollary}
\begin{proof}
By the theorem of Fujita \plaincite[Théorème, p.~473]{Fujita} and Takeuchi \plaincite[Théorème~IV, p.~181]{Takeuchi}, $\Omega$ is Stein and hence admits a complete Kähler metric. The wedge condition says that $g\perp\ker B_{p,p}$, and the lemma gives $B_{p,p}\ge n\operatorname{Id}$ on this orthogonal complement. Hence $M(g)\le n^{-1}\|g\|_{L^2}^2$ in Theorem~\ref*{thm:demailly}.
\end{proof}

\subsection{A Poisson equation on a compact Kähler manifold}

By Buchdahl \plaincite[\S1, p.~289]{Buchdahl}, the operator $u\mapsto *(i\partial\ddbar u\wedge\omega^{n-1})$ is an isomorphism between the mean-zero subspaces of $W^{2,2}$ and $L^2$. Taking $w=u\omega^{p-1}$ gives the following lemma.

\begin{lemma}\label{lem:poisson}
Let $(M,\omega)$ be a compact connected Kähler manifold of dimension $n$, let $1\le p\le n$, and let $a\in L^2(M)$ satisfy $\int_M a\,dV_\omega=0$. There exists a form $w\in W^{2,2}_{p-1,p-1}(M)$ satisfying
\[
 i\partial\ddbar w\wedge\omega^{n-p}=a\,dV_\omega
\]
in the sense of distributions, with $\|w\|_{W^{2,2}}\le C_{M,\omega,p}\|a\|_{L^2}$.
\end{lemma}

\section{Proof of the Main Theorem}

For an open set $D$ in an affine chart, write $A^2(D)=\mathcal O(D)\cap L^2(D,dV)$, with respect to Euclidean volume.
\begin{lemma}\label{lem:berg}
Let $\Omega\subset\PP^n$ be a pseudoconvex domain whose complement $K$ is not locally pluripolar. There exist a ball $B$ relatively compact in an affine chart and a function $h\in A^2(B\cap\Omega)$ that has no holomorphic extension to $B$.
\end{lemma}
\begin{proof}
Choose an affine chart $U\simeq\CC^n$ containing a point at which $K$ is not locally pluripolar, and fix $a\in U\cap\Omega$. The set $U\cap\Omega$ is connected, since it is obtained from $\Omega$ by removing its intersection with a complex hyperplane. Since $U\cap\Omega$ is pseudoconvex and $K\cap U$ is nonpluripolar, Sadullaev's slicing theorem \plaincite[Theorem~1, p.~92]{Sad} gives a complex line $\ell\subset U$ through $a$ such that $K\cap\ell$ is nonpolar. Choose a ball $B\Subset U$ containing $a$ such that $K\cap B\cap\ell$ is nonpolar in the disc $\Delta=B\cap\ell$.

By Björn's removability theorem \plaincite[Theorem~13.3, p.~217]{Bj}, there is an $h_0\in A^2(\Delta\cap\Omega)$ that does not extend holomorphically to $\Delta$. If $\Delta\cap\Omega$ is disconnected, one can simply take the indicator of one of its components. For $n\ge2$, successive applications of the Ohsawa--Takegoshi theorem \begin{NoHyper}\cite{OT}\end{NoHyper} on hyperplane sections of the bounded pseudoconvex open set $B\cap\Omega$ extend $h_0$ to a function $h\in A^2(B\cap\Omega)$. For $n=1$, take $h=h_0$. Any holomorphic extension of $h$ to $B$ would restrict to a holomorphic extension of $h_0$ to $\Delta$.
\end{proof}

Fix $1\le p\le n$ and set
\[
 I_p(f)=\int_\Omega f\wedge\omega^{n-p},\qquad f\in L^2_{p,p}(\Omega).
\]
Since $\Vol(\Omega)<\infty$, $I_p$ is a bounded linear functional.

\begin{lemma}\label{lem:mean}
Suppose that $K=\PP^n\setminus\Omega$ is not locally pluripolar. For every
$f\in L^2_{p,p}(\Omega)$ there exists a correction form $v_{\mathrm{corr}}\in\Dom\ddbar_{p,p-1}$ such that
\begin{equation}\label{eq:correction}
 I_p(\ddbar v_{\mathrm{corr}})=I_p(f).
\end{equation}
\end{lemma}
\begin{proof}
It suffices to find one $\tau\in\Dom\ddbar_{p,p-1}$ with $I_p(\ddbar\tau)\ne0$.
Indeed, $v_{\mathrm{corr}}=I_p(f)\tau/I_p(\ddbar\tau)$ then satisfies~\textup{(\ref*{eq:correction})}.

Choose $B$ and $h$ as in Lemma~\ref*{lem:berg}, and put $D=B\cap\Omega$. The function $h$ has no holomorphic extension to $B$. It is also square-integrable for the Fubini--Study volume, which is comparable to Euclidean volume on $B$.

The zero extension $\widetilde h$ to $B$ therefore satisfies $\ddbar\widetilde h\ne0$ in distributions. Wedging with $\omega^{n-1}$ identifies $(1,0)$-forms with $(n,n-1)$-forms. Testing $\ddbar\widetilde h$ therefore gives a
$\theta\in\mathcal A_c^{1,0}(B)$ such that
\[
 \int_D h\,\ddbar\theta\wedge\omega^{n-1}\ne0 .
\]
Set $\tau=h\theta\wedge\omega^{p-1}$ on $D$ and extend it by zero to
$\Omega$. The compact support of $\theta$ in $B$ and the holomorphicity
of $h$ on $D$ give $\tau\in\Dom\ddbar_{p,p-1}$. The displayed integral
is $I_p(\ddbar\tau)$.
\end{proof}

\begin{proof}[Proof of Theorem~\ref*{thm:main}]
For $\condref{nonpluripolar}\Rightarrow\condref{all}$, fix $1\le p\le n$ and let
$f\in L^2_{p,p}(\Omega)$ be $\ddbar$-closed. By Corollary~\ref*{cor:orth},
it suffices to find $v\in\Dom\ddbar_{p,p-1}$ with
$(f-\ddbar v)\wedge\omega^{n-p}=0$.
Choose $v_{\mathrm{corr}}$ as in Lemma~\ref*{lem:mean}. The zero extension $F$ of
$(f-\ddbar v_{\mathrm{corr}})\wedge\omega^{n-p}$ to $\PP^n$ is an $L^2$ top-degree form
of integral zero. Lemma~\ref*{lem:poisson} gives a $(p-1,p-1)$-form
$w\in W^{2,2}(\PP^n)$ such that $i\partial\ddbar w\wedge\omega^{n-p}=F$.
Thus $v=v_{\mathrm{corr}}-i\partial w|_\Omega$ has the required property.
When $p=n$, we have $\omega^{n-p}=1$, so this gives $f=\ddbar v$ directly.

If $K$ is locally polar, restriction identifies the two $L^2$ complexes isometrically, as shown in the introduction using \plaincite[Corollary~2.4]{CWW}. Consequently,
\[
 \hL pp(\Omega)\simeq H^{p,p}(\PP^n)=\CC[\omega^p],\qquad 1\le p\le n .
\]
This proves $\condref{some}\Rightarrow\condref{nonpolar}$. Since \condref{nonpolar} and \condref{nonpluripolar} are equivalent and $\condref{all}\Rightarrow\condref{some}$ is immediate, the proof is complete.
\end{proof}

\section{Examples}\label{sec:examples}

\begin{example}[A Fatou--Bieberbach domain]
For $n\ge2$, projective space $\PP^n$ contains two domains biholomorphic to $\CC^n$ with different diagonal $L^2$ Dolbeault cohomology for the induced Fubini--Study metric:
\begin{itemize}[label=\textbullet,leftmargin=1.6em]
\item Let $H_\infty$ be the hyperplane at infinity and $U=\PP^n\setminus H_\infty\simeq\CC^n$. Since $H_\infty$ is locally pluripolar, $\hL pp(U)=\CC[\omega^p|_U]$ for $1\le p\le n$.
\item There is a domain $V\subset U\setminus\overline{\mathbb B^n}$ biholomorphic to $\CC^n$, by Rosay and Rudin \plaincite[Theorem~8.5, p.~72]{RR}. A proper subdomain of $\CC^n$ biholomorphic to $\CC^n$ is called a \emph{Fatou--Bieberbach domain}. The domain $V$ is Stein, and its projective complement contains the unit ball. Thus Theorem~\ref*{thm:main} gives $\hL pp(V)=0$ for $1\le p\le n$. The boundary of $V$ in $\PP^n$ cannot be Lipschitz. Otherwise Harrington's theorem \plaincite[Theorem~1.1]{Har} would give a bounded negative plurisubharmonic exhaustion of $V$, contrary to the fact that every negative plurisubharmonic function on $\CC^n$ is constant.
\end{itemize}
\end{example}

\begin{samepage}
\begin{example}[The forward escaping set of a Hénon map]\label{ex:henon}
Let $U=\PP^2\setminus H_\infty\simeq\CC^2$. In affine coordinates, a quadratic Hénon map has the form $F_{a,c}(z,w)=(w,w^2+c-az)$, where $a\in\CC^\times$ and $c\in\CC$. Take $F=F_{2,0}$ and write $F^m$ for its $m$th iterate. Its forward escaping set is
\[
 U_+=\bigl\{x\in U:\lVert F^m(x)\rVert\to\infty \text{ as } m\to\infty\bigr\}.
\]
By Bedford and Smillie \plaincite[Proposition~3.3, p.~75, and Corollary~3.5, p.~76]{BS}, $U_+$ is pseudoconvex. Since $\det_{\CC}DF=2$, Friedland and Milnor \plaincite[Lemma~3.7, p.~78]{FM} give $\Vol(U\setminus U_+)=0$.

The set $U\setminus U_+$ is not locally pluripolar \plaincite[paragraph after Proposition~3.9, p.~78]{BS}. Theorem~\ref*{thm:main} applied to $U_+$ and the preceding example for $U$ give
\[
 \hL 11(U_+)=0,\qquad \hL 11(U)=\CC\,[\omega|_U].
\]
Thus restriction $L^2_{p,q}(U)\to L^2_{p,q}(U_+)$ is an isometric isomorphism in every bidegree, but the $L^2$ Dolbeault cohomology groups of $U$ and $U_+$ in bidegree $(1,1)$ are not isomorphic.
\end{example}
\end{samepage}

\end{document}